\documentclass{ip-journal}
\usepackage{amsmath,epsf,amssymb,latexsym,enumerate,cite,shadow,array,color}

\theoremstyle{definition}

\theoremstyle{remark}

\numberwithin{equation}{section}

\def\E#1#2{\mathbb{E}_{#1}[#2]}

\def\IC{\mathbb{C}}
\def\IN{\mathbb{N}}

\def\IR{\mathbb{R}}
\def\IP{\mathbb{P}}
\def\IK{\mathbb{K}}
\def\cal{\mathcal}

\def\CN {{\cal N}}

\def\CD {{\cal D}}
\def\CF {{\cal F}}

\def\CP {{\cal P }}

\def\CL {{\cal L}}

\def\CO {{\cal O}}
\def\CX {{\cal X}}

\def\CE {{\cal E}}

\def\CH {{\cal H}}

\def\CK{{\cal K}}
\def\CW{{\cal W}}
\def\CY{{\cal Y}}

\newcommand{\eq}[1]{Eq.~(\ref{eq:#1})}

\newcommand{\be}{\begin{equation}}
\newcommand{\ee}{\end{equation}}
\newcommand{\ba}{\begin{array}}
\newcommand{\ea}{\end{array}}
\newcommand{\bea}{\begin{eqnarray}}
\newcommand{\eea}{\end{eqnarray}}

\renewcommand{\Im}{{\rm Im }}
\renewcommand{\Re}{{\rm Re }}

\def\one{{\hbox{ 1\kern-.8mm l}}}

\def\Hom{{\mbox{Hom}\,}}

\def\tr{{\rm tr\,}}

\def\p{\partial}
\def\ba{\bar{a}}

\def\bZ{\bar{Z}}

\def\bj{{\bar{j}}}

\def\bl{{\bar{l}}}

\def\bs{{\bar{s}}}

\def\bJ{{\bar{J}}}

\def\bZ{\bar{Z}}

\begin{document}

\title{ Holomorphic feedforward networks}

\author{Michael R. Douglas}
\address{CMSA, Harvard University, Cambridge 02138}

\email{mdouglas@cmsa.fas.harvard.edu}



\date{May 5, 2021.}



\begin{abstract}
A very popular model in machine learning is the feedforward neural network (FFN).
The FFN can approximate general functions and mitigate the curse of dimensionality.
Here we introduce FFNs which represent sections of holomorphic line bundles on complex manifolds,
and ask some questions about their approximating power.  
We also explain formal similarities between the standard approach to supervised learning and
the problem of finding numerical Ricci flat K\"ahler metrics, which allow carrying some ideas between
the two problems.
\end{abstract}

\maketitle

I'm delighted to contribute to this volume for Bernie Shiffman, who along with Steve Zelditch initiated me into complex geometry as it is actually done by mathematicians.  Our papers on the statistics of flux vacua
\cite{douglas_critical_2004,douglas_critical_2004-2,douglas_critical_2006} 
remain important in string theory.  They were also my jumping-off point for many other works on complex geometry, in particular work on balanced metrics and the Tian-Yau-Zelditch-Lu expansion with Semyon Klevtsov \cite{Douglas_2009,douglas_black_2008}, 
and on numerical methods for finding Ricci-flat metrics on Calabi-Yau manifolds
\cite{Douglas_2008,Douglas_2007,Braun_2008,Braun_2008b,douglas_numerical_2020}.

More recently, I have been studying machine learning and artificial intelligence.
My interest in this goes back to my graduate school days, where I worked with Gerry Sussman
\cite{applegate1986digital}, a pioneer of AI, 
and took classes from John Hopfield on his famous model of neural networks.
Some reminiscences about these times are in \cite{Mumford}.
Although it took many years, machine learning is now at the forefront of applied mathematics
and many other fields.  Its concepts and technologies are being applied to many problems which
at first might not be thought to be related to learning or statistics.

One very active topic is numerical methods for PDE
inspired by machine learning, and a few examples are
\cite{carleo_solving_2017,sirignano_dgm:_2018,han2018solving}.
A good way to explain the relevance of ML
is to recall the ``curse of dimensionality'' \cite{bellman_dynamic_1957}.
A central part of any numerical method for PDE is to define a finite dimensional
space of parameterized functions $\CF$ which are used to approximate the solutions.  
The simplest way to do this is to specify the function values at a fixed subset of points.
This approximation is controlled by restricting the variation of the target function,
for example one can put a uniform bound on its derivatives
such as $|f(x)-f(y)| \le L|x-y|$.  While this gets us down to a finite problem which can be put on a computer,
on a $d$-dimensional manifold $M$ we require $\CO(L^d)$ parameters.  This is a disaster for the
high dimensional PDEs of many-body physics and control theory, and problematic even
for moderate dimensions such as the $d=6$ Calabi-Yau manifolds of interest in string theory.

Machine learning also involves approximating functions, often functions on thousands or even
millions of variables, such as digital images.  Thus the curse of dimensionality is a central problem
and there are many ways to deal with it.
The simplest class of ML problems is supervised
learning, in which one needs to estimate an unknown function $f:M\rightarrow\IR$ from data.  
This includes classification -- say, we take $M$ to be the space of maps from pixels to intensities, 
and we define $f>0$ (resp. $f<0$)
to mean that the image is a picture of a cat (resp. a dog).
What makes this ``learning'' is that we do not say much about $f$ {\it a priori}, rather we tell the
computer how to infer $f$ from examples.  Thus, we give it
a dataset $\CD$ of input-output pairs $(x_i,f(x_i))$ sampled from some underlying distribution $\rho$ on $M$.
We then postulate a class of functions $\CF$ (a model), and an algorithm which looks at the data and
produces a $\hat f\in\CF$ which generalizes to the underlying distribution, {\it i.e.} $\E{\rho}{\hat f-f}$ is
expected 
to be small.   A straightforward
way to try to do this is to numerically minimize the error $\E{\CD}{\hat f-f}$ on the dataset.  

In traditional statistics, one postulates a simple $\CF$, say functions which are a linear combination
of $N$ fixed basis functions.  One then finds a tradeoff between accuracy of approximation (which increases with $N$)
and errors arising from fitting the randomness in sampling from $\rho$ (which also increase with $N$).
This is quantified by concepts such as bias-variance tradeoff, and it leads one to prefer 
simple models, those which use the smallest
possible number of parameters $N$ needed to describe the data.

Much of the recent work rejects the idea that $\CF$ should
be simple.   In ``deep learning,'' one takes $\CF$ to be the class of functions which can be realized by
a class of feedforward neural networks with a definite number of layers, parameters, {\it etc}..
These functions depend nonlinearly on the parameters in fairly complicated ways, but in return
the curse of dimensionality is mitigated.  While the nonlinearity
of $\CF$ would be expected to lead to a bad optimization problem with many local minima, in practice this
turns out not to be the case.   Indeed,
by taking a large number of parameters, one can simplify the optimization landscape.  
Given enough parameters (many more than the number of data points), one can find a $\hat f$ which precisely
fits the data points.  For such a $\hat f$, $\E{\CD}{\hat f-f}=0$ and it turns out to be easy to find such minima.
However, the dogmas of statistics tell us that such minima are ``overfit'' and give a poor
description of the underlying distribution, {\it i.e.} $\E{\rho}{\hat f-f}$ is large.

But surprisingly, with the right choices of $\CF$ and optimization algorithm, this is not the case.
Such models can generalize well,
apparently violating this dogma.  Understanding how this can be is a very active
field of research, of which a few recent works include
\cite{zhang_understanding_2016,belkin2019reconciling,bartlett_deep_2021}.

Carrying over these advantages from ML to PDE, one could hope that a
function class $\CF$ produced by a feedforward neural network would have the same advantages in 
numerical solution of PDE, mitigating the curse of dimensionality and simplifying the optimization problem.
These hopes have been borne out in many works such as 
\cite{carleo_solving_2017,sirignano_dgm:_2018,han2018solving}.

The present work is inspired by work on numerical methods to find Ricci flat metrics on hypersurfaces in
projective space
\cite{donaldson2009some,Douglas_2008,Headrick_2013}.  
Here $\CF$ is a class of K\"ahler metrics, and in previous works this was taken to be the space of
all metrics which could be obtained by pullback from a Fubini-Study metric, thus log linear in the parameters.
We adapt the feedforward neural network to this problem and define ``holomorphic networks''
and ``bihomogeneous networks,'' which represent subsets of these metrics
using a nonlinear parameterization.  In \cite{douglas_numerical_2020}
we study this approach numerically.  Here we point out some formal parallels with machine learning, and
ask some mathematical questions about the ability of these FFNs to
approximate functions.

\section{ Network approximations to K\"ahler metrics }

\subsection{ Feed-forward networks }

A feed-forward network (FFN, also called
MLP for multilayer perceptron) is a nonlinear map $F[{\vec W}]$ from a vector space $\CX$ 
to another vector space $\CY$ with parameters $\vec W$.
It is built by composing a sequence of maps
which alternate between two types, general linear maps $W$ and fixed nonlinear transformations $\theta$, as in
\be \label{eq:defMLP}
F[{\vec W}] = W^{(d)} \circ \theta|_{V_{d-1}}  \circ W^{(d-1)} \circ \ldots \circ
\theta|_{V_2}  \circ W^{(2)} \circ \theta|_{V_1} \circ W^{(1)} .
\ee
Each linear map $W^{(i)}$ has as its range a new vector space $V_i$, so
\bea
W^{(1)} &\in& \Hom( \CX, V_1 ),\\
W^{(2)} &\in& \Hom( V_1, V_2 ), \\
\vdots \nonumber\\
W^{(d)} &\in& \Hom( V_{d-1}, \CY )
\eea
The combination $\theta\circ W$ is called a layer, with 
the final layer $W^{(d)}$ being an exception in not having $\theta$.
The number of layers $d$ is the depth.

Usually one takes the domain $\CX \cong \IK^D$ and range $\CY\cong \IK^{D'}$ with $\IK\cong\IR$,
though $\IK\cong\IC$ is sometimes used in signal processing applications.
We are free to choose the dimensions $D_i$ of the intermediate spaces $V_i$, 
the width hyperparameters of the network.\footnote{ The term hyperparameter in machine learning
refers to a choice which is not determined automatically by optimization, in contrast to parameters
(also called ``coefficients''
or ``weights'') which are generally chosen to minimize an energy (``objective'' or ``loss'') function.}
By ``the width'' we will mean $\max D_i$.

Generally one allows the $W$'s to be arbitrary linear transformations, so the parameters consist of
the list of $W^{(i)}$.
Thus, defining $D=\dim \CX$ and $D'=\dim \CY$, we have
\be \label{eq:defCW}
\vec W \in \CW_{\vec D} \cong \oplus_{i=1}^{d} \Hom( \IK^{D_{i-1}} , \IK^{D_{i}} ) .
\ee

To define the $\theta$'s, we start with the one dimensional case $\theta|_\IK:\IK\rightarrow\IK$,
which is called the activation function.  This could be any function;
two popular choices for $\IK=\IR$ are
$\theta(x)=\tanh x$, and the ``rectified linear unit'' or ReLU function 
\be \label{eq:ReLU}
\theta_{ReLU}(x) = \begin{cases} x, x\ge 0 \\ 0, x<0 \end{cases}.
\ee
To define $\theta_V$ for a general vector space $V$, we pick a basis $e_i$ for $V$ and apply
$\theta|_\IK$ componentwise,
\be \label{eq:thetabasis}
\theta_V\left( \sum_i c_i\; e_i \right) = \sum_i \theta_\IK(c_i)\, e_i.
\ee
While this depends on the choice of basis, since every $\theta$ in \eq{defMLP} appears
both prefixed and postfixed by a general linear transformation, 
the parameterized space of maps is independent of this choice.
Thus, for each choice of specific dimensions $\vec D$ and activation function $\theta|_\IK$, we
get a space of FFN maps of dimension $DD_1+D_1D_2+\ldots+D_{d-1}D'$,
which we denote $\CF_{\vec D;\theta}$.

This definition was inspired by a simple model in neuroscience (the ``perceptron''),
in which a neuron is modeled by the computation which
takes an input vector $v$ and produces a single component of $(\theta\circ W)(v)$.
These days neuroscience uses far more sophisticated models of neurons, but 
FFN's are widely used in applied math and machine learning
as function approximating spaces, analogous to spaces
of polynomials or other linear combinations of basis functions.
Thus rather than ``neuron'' one usually uses the term ``unit,''
and the network $\CF[\vec D;\theta]$ contains $D_1+\ldots+D_{d-1}+D'$ units.

It has been shown that feed-forward networks can approximate arbitrary real valued functions.
This is the case even for $d=2$ \cite{cybenko1989approximation}, 
but in this case one can need an exponentially
large number of units, as would be the case for simpler methods of interpolation (the ``curse of dimensionality'').
By using more layers, one can gain many advantages -- complicated functions can be represented with
many fewer units, and local optimization techniques are much more effective.  
How exactly this works is not
well understood theoretically and there are many interesting observations and hypotheses as to how these
advantages arise.

\subsection{ Multilayer holomorphic embeddings }
\label{ss:mhe}

Our goal is to use the FFN \eq{defMLP} to define a parameterized space of functions from a 
projective manifold $M$ to a space $\CY$.  To begin, 
let $\CL\rightarrow M$ be a holomorphic line bundle,
let $\underline{s}=(s_0,\ldots,s_N)$ be a basis of sections of $H^0(\CL)$,
and let $\iota_s$ be the corresponding Kodaira embedding of $M$ into $\IC\IP^N$.

To use \eq{defMLP}, we want to regard this embedding as a map into $\CX\cong \IC^{N+1}$.  We can do this
in a patch $U\subset M$ by choosing a local frame. 
If the maps $W$ and $\theta$ each have a nice geometric interpretation, the choice of frame should
drop out at the end.  Taking the $V_i$'s to be complex vector spaces,
this is evident for $W$ as each $W\in \Hom(V,V')$
corresponds to a holomorphic map $\IP V\rightarrow \IP V'$.

An arbitrary choice of activation function $\theta$ will not have a simple geometric interpretation.
To find one which does, consider the particular case of $M=\IC\IP^n$ and $\CL=\CO(k)$, so that
the sections $s$ are degree $k$ polynomials.  Clearly we want $\theta(s)$ to also be a homogeneous 
polynomial, so the natural choice is
\be\label{eq:needtheta}
\theta( s ) \equiv s^p
\ee
and its componentwise analog \eq{thetabasis}.
Taking into account the dependence on the frame, the geometric interpretation of this is a map
\be
\theta : \underbrace{H^0(\CL) \oplus \ldots \oplus H^0(\CL)}_\text{D times} \rightarrow 
\underbrace{H^0(\CL^p) \oplus \ldots \oplus H^0(\CL^p)}_\text{D times} .
\ee

Another way to relate \eq{defMLP} to geometry is to appeal to the ``metatheorem for vector bundles''
(\cite{huybrechts2005complex}, Thm. 2.2.3), that any canonical construction in linear algebra
gives rise to a geometric version for holomorphic vector bundles.
To use this we define $\theta$ in terms of a canonical multilinear map on $\IK$,
namely the tensor product
\be \label{eq:tensor}
\otimes_p : \IK^m \rightarrow \IK : (z_1,z_2,\ldots,z_k) \rightarrow z_1\;z_2\;\ldots\;z_p .
\ee
By evaluating this on $p$ copies of a section, we get a natural map from the line bundle $\CL^m$ 
to its $p$'th power $\CL^{pm}$.

This interpretation of \eq{defMLP} was made in \cite{douglas_numerical_2020} in the following concrete form.
We took $M$ to be a hypersurface in $\IC\IP^n$, and $\CL=\CO(1)$.  
We then took as our computational representation
of $M$ the image in $H^0(\CO(1))$ of a set of randomly chosen points in $M$.
We could then use the implementations of the $W$ and $\theta$ maps in a standard ML package (Tensorflow). 

While we will generally use $p=2$ in our networks,
one could take a different $p$ in each layer, so the network \eq{defMLP} depends on a vector 
$\vec p=(p_1,\ldots, p_d)$.
As a further generalization,
instead of identifying the input sections $z_i$ in \eq{tensor}, one could multiply the outputs of distinct FFN's,
or make other geometrically consistent combinations.  
Each possible ``architecture'' is naturally associated to a directed graph, in which 
a particular use of \eq{tensor} corresponds to a vertex, the inputs $z_i$ correspond to incoming edges,
and the output $z_1\;\ldots\;z_p$ corresponds to a single outgoing edge.  These are 
examples of tensor networks, about which we will say more below.


The upshot is a construction of parameterized maps from
$M$ into $\oplus_D H^0(\CL^m)$, 
defined by \eq{defMLP} and \eq{needtheta} with parameter space \eq{defCW}.
Let us denote these maps 
as $F^{h}_{\CL;\vec p;\vec D}[{\vec W}]$ with $\vec W\in\CW_{\vec D}$ (or just $F^h[{\vec W}]$ in context).
We can regard them as defining embeddings of $M$ into projective space, which
for sufficiently large $D$ include all possible embeddings by a complete space of sections.

\subsection{ Algebraic metrics and holomorphic networks }

Next, given an embedding of $M$ into projective space, we can pull back the Fubini-Study metric
to get a K\"ahler metric on $M$. 
Let $\CK^{FS}_{\CL}$ be the space of Fubini-Study metrics on $\IC\IP^N\cong \IP H^0(\CL)^*$, parameterized
by an $(N+1)\times (N+1)$ positive definite hermitian matrix $G_{I,\bJ}$.  The K\"ahler potentials
\be \label{eq:a1}
K^{FS}[G] = \log \sum_{I,\bJ} G_{I,\bJ} s^I \bs^\bJ .
\ee
then define K\"ahler metrics on $M$, the space of which we also denote
$\CK^{FS}_{\CL}$.  Replacing $\CL$ with
the series of line bundles $\CL^{k}$, we get embeddings into higher
dimensional projective spaces, and a sequence of spaces of Fubini-Study metrics $\CK^{FS}_{\CL^k}$.
This gives us the ability to approximate a given K\"ahler metric to arbitrary accuracy, by taking sufficiently high $k$.
As discussed in \cite{tian1990set,donaldson2009some}, this accuracy can go as $k^{-\nu}$ for any $\nu$.

This is a very nice parameterized class of metrics, and we will shortly discuss a method for finding
the best numerical approximation to the Ricci flat metric in this class.
However, there is a problem with this numerical application.
Since $N\sim k^{\dim_\IC M}$, the number of coefficients will go as $k^{\dim_\IR M}$, 
the curse of dimensionality referred to earlier.

Let us try to break this curse by replacing the complete space of sections with the image of 
$F^h_{\CL;\vec p;\vec D}[{\vec W}]$, to produce a space of metrics $\CK^{h}_{\CL;\vec p;\vec D}$.
Now the parameter spaces \eq{defCW} for the construction in \S \ref{ss:mhe} come in a
variety of sizes, obtained by adjusting the intermediate dimensions $\{D_i\}$.  And since the degree
$k$ is related to the depth as $k=\prod_{i=1}^{d-1} p_i$ (in terms of $\vec p$ in \eq{needtheta}), 
the number of coefficients scales with $k$ only as $\CO(D^2\log k)$.

We now write this out explicitly for $d=2$ and $p=2$.  A slightly 
simpler expression is obtained by removing the final weight matrix $W^{(2)}$
from \eq{defMLP}, which is redundant if we can vary $G$ in \eq{a1}.  We then have
\be\label{eq:holo2net}
K^{h}[W^{(1)},G] = \log \sum_{1\le I,\bJ \le D_1} G_{I,\bJ} (\sum_i W^{(1),I}_i s^i)^2 (\sum_j W^{(1)*,\bJ}_\bj \bs^\bj)^2 .
\ee
The nonlinear dependence on the weights $W$ is a general feature of FFN's.

We studied this construction computationally along lines we describe below,
and found that it is not very good at approximating Ricci flat metrics.
This may be because these metrics are pullbacks of degenerate FS metrics on the embedding space,
called partial Fubini-Study metrics (see for example \cite{zelditch2013pointwise}).

\subsection{ Bihomogeneous networks }
\label{ss:bihom}

A variation on the previous construction is to first make sesquilinear combinations of the holomorphic
sections and then apply the nonlinear network construction.  Thus, we could take
\be
\CX \equiv H^0(\CL) \otimes \bar{H^0(\CL)} ,
\ee
embed $M\hookrightarrow \CX$ by taking the outer product $s^I\times \bs^\bJ$ of the basis of sections
with its complex conjugate,
and use the network \eq{defMLP} to define parameterized functions $F[\vec W]:\CX\rightarrow\IR$.
We then replace the K\"ahler potential \eq{a1} with
\bea\label{eq:Kb}
&K^{b}_{\CL;\vec p;\vec D}[ \vec W] = \log F_{b}  \\
\label{eq:Fb}
& F^{b} =  W^{(d)} \circ \theta  \circ W^{(d-1)} \circ \ldots \circ
 \theta \circ W^{(1)} \cdot (\Re\; s^I\bs^\bJ, \Im\; s^I\bs^\bJ) .
\eea
In terms of a concrete embedding by polynomials, this amounts to taking the linear spaces $V_i$
to be subspaces of the space of bihomogeneous polynomials of degree $(n_i,n_i)$, while
$\theta$ is still \eq{needtheta}. 

We will denote a particular K\"ahler potential of this type as $K^b[\vec W]$ and the space of these as
$\CK^{b}_{\CL,\vec p,\vec D}$.  Here $D=D_0=b^0(\CL)^2$ and $D'=D_d=1$, with the intermediate $D_i$ adjustable.
Note that $\CX$ has a real structure, and we can take the intermediate values and weights in \eq{defMLP}
to be real.  Thus this network has $DD_1+\ldots+D_{d-1}$ real parameters.
For sufficiently large $D_i$, this reproduces the complete space $\CK^{FS}_{\CL^{p_1\ldots p_{d-1}}}$
of FS  K\"ahler potentials, though with a more complicated parameterization. 

Here is a depth 2 bihomogeneous network with $p=2$ (the analog of \eq{holo2net}), and a depth 3 network:
\bea\label{eq:bihom2net}
K^b_{\CL;2;D_1} &=& \log \sum_{1\le I \le D_1} W^{(2)}_{I} (\sum_i W^{(1),I}_{i,\bj} s^i \bs^\bj)^2 \\
\label{eq:bihom3net}
\qquad K^b_{\CL;\vec 2;D_1,D_2} &=& \log \sum_{1\le I \le D_2} W^{(3)}_{I} \left( \sum_{1\le J \le D_1} W^{(2),I}_{J}
(\sum_i W^{(1),J}_{i,\bj} s^i \bs^\bj)^2 \right)^2 .
\eea

These networks are better behaved than the holomorphic networks for small intermediate widths.
In particular, the function class with $d$ layers is simply
contained in the function class with $d+1$ layers.  To be precise, take $p=2$, then
we apply a rescaling (denoted $2\times$) to match the K\"ahler classes from the two constructions,
\be
2\times\CK[\CL;\vec 2;\vec D] \subset \CK[\CL;\vec 2;\vec D \oplus D_{d}] .
\ee  
For $\CK^h$ one needs $D_{d}\ge h^0(\CL^{2^{d-1}})$.
But for $\CK^b$, this can be accomplished with $D_{d}=1$ (for example, $D_2=1$ in \eq{bihom3net}).

In particular one can start with the 1 layer network $F= |W \cdot z|^2$ and repeatedly apply $\theta \cdot {\bf id}$ to get 
$F= |W \cdot z|^{2k}$, a function sharply peaked at $z\propto W$.  One can then add this to another
network $F'$, increasing its widths by one and adding a feature with size $1/k$.  Thus one can describe
structure on the same short length scales as \eq{a1} but using many fewer parameters.

In \cite{douglas_numerical_2020} we give numerical evidence
that these metrics are quite good at representing the Ricci flat metrics.
Again, there are variations of this construction labeled by directed graphs.
We could also take differences of these potentials to define classes of real valued functions.

\subsection{ Matrix product states }
\label{ss:mps}

These were originally proposed to describe spin chains in quantum physics,
and have been applied to ML problems in \cite{anandkumar2014tensor,stokes2019probabilistic}
and other works.  
We will use them to make a third definition of a parameterized subset of the Fubini-Study metrics.

Without going into all the details, a spin is a quantum system whose
wave function is a point in $\IP V^*$ with $V\cong\IC^D$, while a chain of $N$ spins has a wave function $\Psi$ in 
$\IP\CX$ where $\CX\equiv\otimes^N V^*$. 
Since $N$ might be the number of atoms in a macroscopic object, the curse of dimensionality is a dominant
aspect of these problems.  Thus the need for parameterized low dimensional subsets of $\CX$ is even greater.

A matrix product state (MPS) is defined by specifying a series of auxiliary linear spaces $L_i$,
and a series $T_i$ of multilinear maps
\be
T_i : V \rightarrow \mbox{Hom}(L_i,L_{i+1}) .
\ee
Typically one takes all of the $L_i$ to be isomorphic, and then writes
\be
\Psi_T = \mbox{Tr}\; T_1 T_2\ldots T_N .
\ee
This gives us an $N(\dim L)^2-1$ dimensional subset of the full $D^N-1$ dimensional space of wave functions.

To adapt this to the problem at hand, we could take $V\equiv H^0(\CL)$ to get a parameterized
subset of holomorphic sections of $\mbox{Sym}^N V$.  
While we have not yet tried
this construction or its many variations, one which seems particularly natural is a bihomogeneous
version with
\be
K^T \equiv \log \mbox{Tr}\; (T_1\cdot z) (T_1\cdot z)^\dag (T_2\cdot z) (T_2\cdot z)^\dag \ldots
(T_d\cdot z) (T_d\cdot z)^\dag .
\ee

\subsection{ Numerical methods }
\label{ss:num}

Suppose we want to approximate a given K\"ahler potential $K$ with one from these classes, say $\CK^b$.
This is a problem in interpolation of functions, which is closely related to supervised machine learning.  
To get a dataset, we sample points $x_i$ from $\CX$ and evaluate $y_i=K(x_i)$, to get ``data'' $(x_i,y_i)$.
We then postulate a ``loss function,''
a measure of the distance between the dataset and the model $(x_i,K_b(x_i))$.
This is a function of the parameters of $K^b$, which we then minimize using computational methods,
usually gradient descent.

Now, the problem of finding a Ricci flat K\"ahler metric can be phrased as an approximation problem.
Since the first Chern class vanishes, the canonical bundle has a unique section up to overall scale.
This is a nonvanishing holomorphic $d$-form, call it $\Omega$.  Then the Ricci flat K\"ahler form $\omega$
satisfies
\be\label{eq:rflat}
\omega^d = \CN\Omega\wedge\bar\Omega
\ee
where the constant $\CN$ can be determined by integrating both sides.  One could then use
the $L^p$ norm of the difference as a loss function, as in \cite{Headrick_2013}.

In \cite{douglas_numerical_2020} we carried this out to find numerical Ricci flat metrics on
quintic hypersurfaces of dimension $3$, following the general approach developed in 
\cite{Douglas_2008,Headrick_2013}.  To sample from a hypersurface in $\IC\IP^{d+1}$, we
sample pairs of points using some measure $\mu$ on $\IC\IP^{d+1}$.
We then find the intersection of the corresponding line
with the hypersurface.  By a result of Shiffman and Zelditch (\cite{shiffman_distribution_1999}, lemma 3.1),
the resulting set of points is distributed by the pullback of $\mu$.
Then, since one can compute the right hand side of \eq{rflat}, 
one can treat this as a problem in function interpolation,
choosing the weights in \eq{Kb} to best fit this given function.  This is a task well suited to machine
learning software.

\section{ Accuracy of approximation }

What determines the accuracy of the
FFN description \eq{Kb} of a canonical K\"ahler metric?  Is it the depth of the network or the number of its parameters?
Here is a precise form of the question and some speculations.

\subsection{ Approximation by algebraic metrics }
To begin, recall the discussion of approximation of metrics from \cite{tian1990set,donaldson2009some}.
Given any K\"ahler metric $\omega$ on $M$, we can find a sequence of algebraic metrics
(pullbacks of degree $k$ Fubini-Study metrics \eq{a1} to $M$) such that
\be \label{eq:approx}
||\omega-\omega_k|| = o(k^{-\nu})
\ee
for any $\nu$.  
This can be shown by using the Tian-Yau-Zelditch-Lu expansion
for the density of states on $M$,
\bea\label{eq:defrho}
\rho_{h}(z) &\equiv& \sum_I |s^I(z)|^2_h \\
&=& 1 + k^{-1} a_1(\omega) + k^{-2} a_2(\omega) + \ldots , \label{eq:TYZL}
\eea
where $h$ is a Hermitian metric on $L^k$ with curvature $-2\pi ik\omega_h$,
and the $a_i$ are local invariants constructed from $\omega$.
To define $\rho$ we introduce a Hermitian metric on $H^0(L^k)$,
\be\label{eq:hilb}
||s||^2_{Hilb(h)} = R \int_M d\mu_h \,|s|_h^2 ,
\ee
where $d\mu_h=\omega_h^n/n!$.
The sum in \eq{TYZL} is taken over an orthonormal basis with respect to \eq{hilb}.
Thus it depends on $h$ both through the inner product \eq{hilb} and explicitly in $|\ldots|_h$.

Now, starting with $h$ such that $\omega_h=\omega$, we can find a $\tilde h$ whose curvature is
\be\label{eq:astep}
\omega_{\tilde h} = \omega_h + k^{-1} i \partial\bar\partial\log(\rho_h).
\ee
One way to see this is to regard the sum in \eq{defrho}
as a bihomogeneous polynomial in the sense of \S \ref{ss:bihom}.  As such it can be
directly interpreted as $\tilde h$.  Taking $\partial\bar\partial\log$ of \eq{TYZL} then gives us \eq{astep}.

We can now define the balanced metrics.  These have $\tilde h=h$ in \eq{astep}, and thus
\be\label{eq:baldef}
\rho_{\tilde h}(z) = 1.
\ee
When a balanced metric exists, 
by \eq{TYZL} it is a canonical approximation to a Ricci flat or constant scalar curvature metric.

The ideas above can be combined to show \eq{approx}.  One can use \eq{TYZL} to compute the leading term
in $\rho_{\tilde h}$ in an expansion in $1/k$, and then \eq{astep} to define a corrected $\tilde h'$ for which
this term agrees with the one in $\rho_h$.  This procedure can be carried out to arbitrary order.
The problem is that one needs $k^{2n}$ coefficients to describe this sequence.
This fits with the expectation from ``curse of dimensionality'' for
a basis of functions which is localized on a length scale $L\sim 1/k$.

\subsection{ Approximation by bihomogeneous metrics }
\label{ss:netmet}

Could \eq{approx} hold 
with $\omega_k$ taken from some restricted class of metrics?
In particular, consider the subsets of metrics $\CK^{b}_{\CL;\vec p;\vec D}$ defined by 
the bihomogeneous networks \eq{Kb}.
For sufficiently large widths $\vec D$, these include the general Fubini-Study metric, and as
we decrease the widths we get subsets which use many fewer parameters,
of order width $\max D_i$ times depth $d$.
Furthermore, as we explained these metrics can contain structure
on the same length scale $L\sim 1/k$ as the general Fubini-Study metric.
Could these networks break the curse of dimensionality?

To sharpen the question,
define a polynomial width sequence of networks to be a sequence of networks with $d\in\IN$,
such that the width $\max D_i$ grows at most polynomially in $d$.  Take $p=2$ for definiteness, then
these networks define sequences of metrics
$\omega_k$ with $k=2^{d-1}$.  Define $\omega$ to be of polynomial complexity if a sequence
exists satisfying \eq{approx} for some $\nu>0$.
We would like to know whether there exist metrics $\omega$ which 
are not of polynomial complexity.

Let us try to bound the maximum distance from the algebraic metrics to an approximate metric, as
\be \label{eq:netbound}
\forall K\in \CK_{\CL^k} , \exists {\tilde K\in \CK^{b}_{\CL;\vec p;\vec D}} \,{\it s.t.}\, 
||\omega_K- \omega_{\tilde K}|| \le B_{\CL^k;\vec p;\vec D}.
\ee
For concreteness we can take $M=\IC\IP^n$, but our arguments here will only be qualitative.
Let us consider a two layer network with $p=2$, so the inputs are the complete set of
degree $(k/2,k/2)$ polynomials, and $\theta(z)=z^2$.  We then want to understand the
dependence $B_{\CL^k;p=2;D_1}$ of the upper bound. 

The simplest hypothesis is that $B=0$ once we match the counting of parameters between $\CK_{\CL^k}$
and $\CK^{b}_{\CL^{k/2};2;D_1}$, in other words when the width $D_1$ satisfies
\be\label{eq:asympDc}
\binom{k+n}{k}^2 = D_1 \binom{\frac{k}{2}+n}{\frac{k}{2}}^2 .
\ee
This has the asymptotic behaviors
\be \label{eq:asympD}
D_1 \sim \begin{cases}
2^{2n} & \text{for } n\ll k, \\
\left(\frac{n}{k}\right)^{k} & \text{for } n \gg k
\end{cases}.
\ee
The intuition is that while we need to combine general quadratic functions of the inputs,
these are only independent for sufficiently high dimension $n$; in low dimension they are linearly dependent.

The qualitative behavior \eq{asympD} can be checked by simplifying even more, to take the inputs to be vectors $X^I$
in $\IR^D$.  For $n=1$, it follows because one can write any degree $4m$ real univariate polynomial as a
difference of squares of two such polynomials.\footnote{ Write $p=g^2-h^2=(g-h)(g+h)$ and then
factor $p$ into terms of degree $1$ and $2$.}  For the $n\gg k$ limit, let $X^I$ be a basis for the
degree $k/2$ polynomials of dimension $D$, then we are trying to minimize an error of the form 
\be \label{eq:fit2}
\CE = || \sum_{I,J} f_{IJ} X^I X^J - \sum_{\alpha=1}^{D_1} (\sum_I W^{(1)}_{\alpha,I} X^I)^2 || .
\ee
If we use the Frobenius norm, this is minimized by regarding $f_{IJ}$ as a matrix and diagonalizing it.
The optimal rank $D_1$ solution is to write $f=c^t\;w\;c$ and keep the top $K$ 
(in magnitude) values of $w$.  The worst case for approximation is all eigenvalues equal,
in which case we need $D_1=D$ reproducing the second case in \eq{asympD}.  In this case one also sees that
taking  $D'<D$ inputs has the same effect.

In the $n\gg k$ limit we even see that the best rank $D_1$ approximation to a quadratic on $\IR^D$ has worst case relative $L_2$ error  $1-D_1/D$.  If we boldly extrapolate this to general $B_{n,k;D_1}$ by replacing $D$ with the critical value
of $D_1$ satisfying \eq{asympDc}, we would conjecture that
\be\label{eq:estB}
B_{n,k;2;D_1} = 
\max 0,\, B_{n,k;1} 
\left(1 - \frac{D_1-1}{D_{match}-1}\right) ;\qquad
D_{match} \equiv  \frac{ \binom{k+n}{k}^2 }{ D_0 } .
\ee
where $D_0=\binom{k/2+n}{k/2}^2 $
and $B_{n,k;1}$ is the worst case error bound for approximating a degree $(k,k)$
polynomial with a single degree $(k/2,k/2)$ polynomial.
In other words, the improvement in approximation ability for a two layer network
is linear in the ratio of the number of available parameters to the minimal
number of parameters for a universal network, meaning one which can exactly reproduce any degree
$k$ polynomial.\footnote{The linearity corresponds to the worst case in which $f$ in \eq{fit2} is the identity.
Constraints on the spectrum of $f$ could lead to better bounds.}

A strong form of the conjecture is that this is even the case if we allow the
choice of inputs to depend on the function we are approximating.
Thus, if we take the inputs to be a subspace of dimension $D_0$ but which can depend on $K$
in \eq{netbound}, the worst case error would still be \eq{estB}.
This is not at all obvious, but the idea is motivated by the effect of taking $D'<D$ in \eq{fit2}.  

This strong form of the conjecture could be generalized to a multilayer network, and suggests that in low dimensions $n$,
the growth $k^{2n}$ is inevitable.  Let us consider a $d$ layer network.
By the above, the final layer can approximate a general function in terms of $2^{2n}$
outputs of independent $d-1$-layer networks.  For example, we can modify \eq{bihom3net} 
by replacing the layer 1 weights $W^{(1),J}_{i,\bj}$ with independent weights $W^{(1),IJ}_{i,\bj}$
for each choice of subnetwork 
indexed by $I$, to get
\be
\label{eq:bihom3tnet}
\qquad K^b_{\CL;\vec 2;D_1,D_2} = \log \sum_{1\le I \le D_2} W^{(3)}_{I} \left( \sum_{1\le J \le D_1} W^{(2),I}_{J}
(\sum_i W^{(1),IJ}_{i,\bj} s^i \bs^\bj)^2 \right)^2 .
\ee
We can repeat this process for each successive layer, producing a
tree structured network (an idea introduced in \cite{e_banach_2020})
whose total number of parameters would be $N \sim 2^{2dn} \sim k^{2n}$.
Since we did not change the assumptions leading to \eq{asympDc}, 
the asymptotic parameter count is still the same as that for a general polynomial.  

Thus the question is, can we improve on this by sharing intermediate results.  
Since the intermediate width of the
tree structured network at layer $l$ is $2^{2n(d-l)}$, in principle there is a lot of scope for this.  
And in the early
layers, the number of independent polynomials is $N \sim (2^l)^{2n}\sim 2^{2nl} \ll 2^{2n(d-l)}$, so for $l<d/2$ one 
might save on parameters by replacing the many subnetworks with a single network which constructs
a complete basis.  However without further assumptions this does not buy much.
Rather than write the analog of \eq{bihom3tnet}, let us just give the modification at layer $d/2$.
Replacing the subnetworks with the general sections of degree $\sqrt{k}$, we take
\be
W^{(d/2),{I_d}{I_{d-1}}\ldots{I_{d/2}}}_{I_{d/2-1}} \Rightarrow
W^{(d/2),{I_d}{I_{d-1}}\ldots{I_{d/2}}}_{i_1\ldots i_{d/2-1},\bj_1\ldots\bj_{d/2-1} } .
\ee
The $2^{dn}$ subnetworks each have $\CO(2^{dn})$ parameters, so the count is still $\CO(k^{2n})$.

Now, it still might be that the subnetworks do not need all of these inputs and one could use fewer parameters.
To study this systematically we would want generalizations of \eq{netbound} to bound the
error for jointly approximating $D_i>1$ functions, given a subspace of the possible inputs.
This second dependence (on the subspace) is crucial as we are trying to drastically reduce the number of inputs 
to the intermediate layers.  But according to the strong form of \eq{estB}, the error would not depend on this
choice; it would still be governed by the ratio of the number of parameters to the number for a universal network.
Any error bound which decreases with $k$ would require a number of parameters which grows as the same power of $k$.

Thus, assuming a very strong conjecture,
we have argued that to obtain \eq{approx} for a general metric 
we need a number of coefficients growing as $k^{2n}$.  
Whether or not we believe this argument, still this would not be too surprising in the worst case.
Could the situation be better in special cases, such as the Ricci flat
and balanced metrics?

Is there a useful definition of a balanced bihomogeneous metric?
Since these are a subset of the embedding metrics, one could
try to define a balanced bihomogeneous metric by finding a functional on embeddings whose minimum
is the balanced metric.  One could then take its minimum over the subset of bihomogeneous metrics.

Indeed, such a functional exists \cite{donaldson2009some},
\be \label{eq:balanced-obj}
\psi(G) = \int_M \nu(z) \,K_G(s(z)) - \frac{1}{N+1}\tr\log G ,
\ee
where $\nu(z)$ is a measure on $M$, and $K$ and $G$ are as in \eq{a1}.
To get the balanced metric defined by \eq{baldef} one must take $\nu$ to be the volume form for this metric.
But the definition makes sense for a general $\nu$, and on a
Calabi-Yau manifold we can use the other canonical volume form
(the r.h.s. of \eq{rflat}) to get a simpler definition which also converges to the Ricci flat metric.

The problem with this is that, as a function of the weights $W$, \eq{balanced-obj} is probably nonconvex.
If one substitutes in \eq{Kb} to get
 an explicit expression for it (analogous to  \eq{fit2}), this will be degree $2^{d-i}$ in the
weights $W^{(i)}$.  Thus it probably does not have a unique minimum.  This problem is familiar from
the neural network literature, but surprisingly enough turns out not to be a serious problem in
the applications to machine learning.  Still, for present purposes an unambiguous definition would be better.
Perhaps there are more geometric definitions of holomorphic network which avoid this problem.

\smallskip
We thank Steve Zelditch 
for suggestions which helped us improve the manuscript.

\bibliographystyle{amsplain}
\bibliography{shiffman5}

\end{document}